\documentclass[conference,compsoc]{IEEEtran}

\ifCLASSOPTIONcompsoc
  \usepackage[nocompress]{cite}
\else
  \usepackage{cite}
\fi
\usepackage{graphicx}
\usepackage{dcolumn}
\usepackage{bm}

\usepackage[utf8]{inputenc}
\usepackage[T1]{fontenc}
\usepackage{mathptmx}
\usepackage{etoolbox}
\usepackage{hyperref}
\usepackage{listings}
\usepackage{amsmath}
\usepackage{amssymb}
\usepackage{tabularx}
\usepackage{makecell}

\ifCLASSINFOpdf
\else
\fi
\makeatother
\begin{document}

\title{Hyperdimensional Computing Provides\\Computational Paradigms for Oscillatory Systems}

\author{\IEEEauthorblockN{Wilkie Olin-Ammentorp}
\IEEEauthorblockA{Mathematics \& Computer Science Division\\
Argonne National Laboratory\\
Lemont, Illinois 60439\\
Email: wolinammentorp@anl.gov}}

\maketitle

\begin{abstract}
  The increasing difficulty in continued development of digital electronic logic has led to a renewed interest in alternative approaches. Oscillatory computing is one such approach that leverages alternative physical systems and computation strategies, but it lacks high-level paradigms for system design and programming. We address this gap by describing a model based on hyperdimensional computing that serves as an ``instruction set'' to integrate oscillatory networks into algorithms for real-valued computing. The expressiveness and compositionality of these instructions allow oscillatory systems to implement both common tasks and novel functions, providing a clear computational role for many emerging hardware devices. We detail the computational primitives of this system, prove how they can be executed via oscillatory systems, quantify the performance of these operations, and apply them to execute multiple tasks including compression, factorization, and classification. 
\end{abstract}

\maketitle

\section{\label{sec:introduction}Introduction}
The rate of progress in computing based on digital electronic logic has slowed, creating a challenge for an information economy accustomed to a rapid pace of hardware improvement \cite{Thompson_2017}. This challenge has motivated numerous responses, including increased hardware specialization, exploration of alternative representations of information, and the development of novel hardware devices \cite{Hennessy_Patterson_2019,Dally_2015,IRDS2022}. 

\begin{figure}[htbp]
    \centerline{\includegraphics[width=0.5\textwidth]{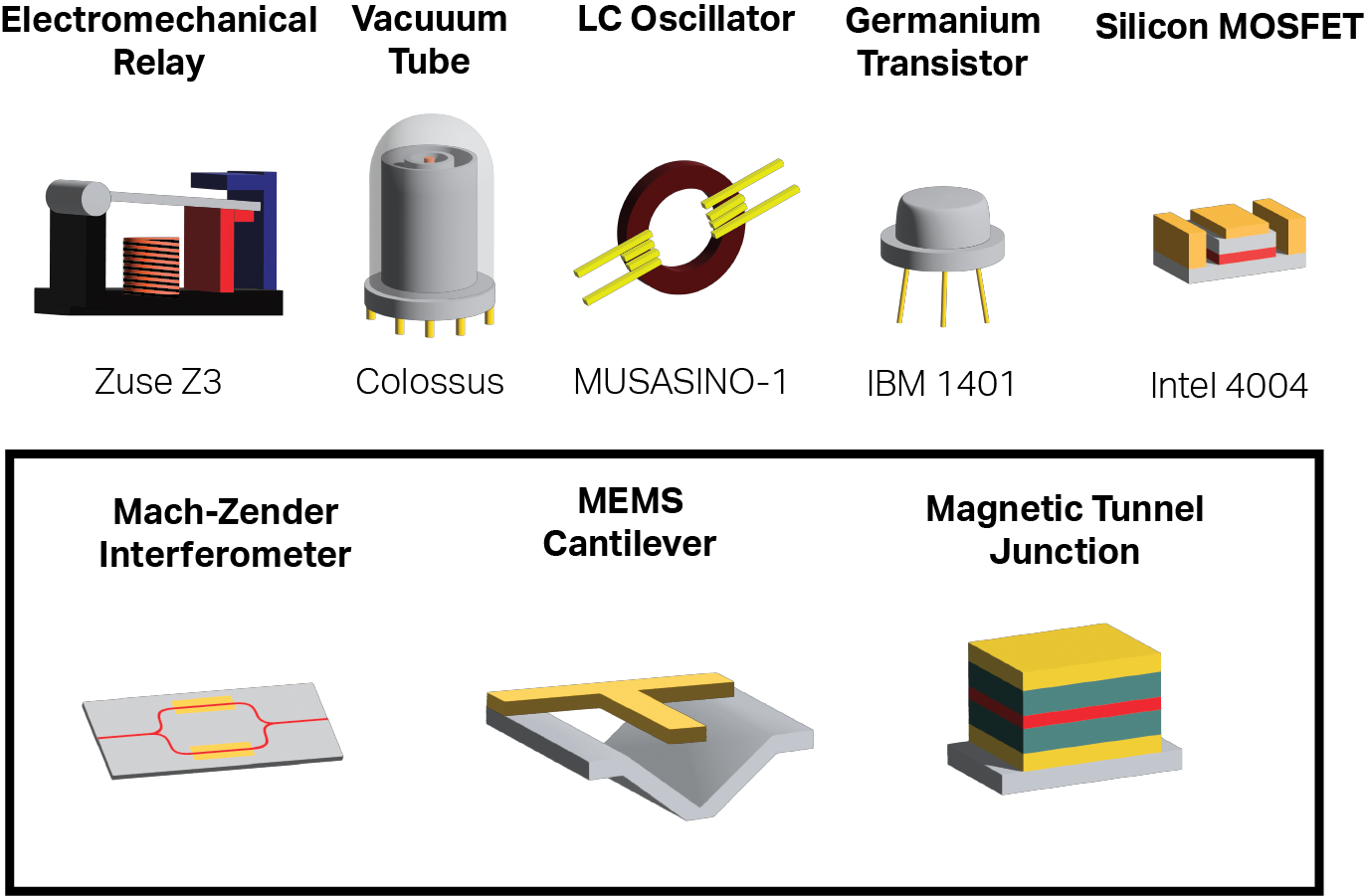}}
    \caption{Illustration of artificial devices employed to represent abstract values. Many devices (top text) have been employed to construct practical computing devices (bottom label), but other emerging devices are still being investigated for this use. We consider a model of real-valued computation employing phase angles. These angles can be expressed by devices that oscillate in the electromagnetic, mechanical, or other domains (inner box).}
    \label{fig1}
\end{figure}

Among these approaches, the development of novel hardware devices offers unique opportunities to utilize and manipulate alternative physical representations of information. These include physical properties such as optical polarization and electronic spin (Figure \ref{fig1}) \cite{IRDS2022}. Utilizing these alternative physical properties can yield systems with entirely different scaling properties from those of traditional digital electronics. For instance, photonic systems demonstrate the ability to represent multiple pieces of information in the same physical device  via wavelength-division multiplexing as well as the ability to carry out multiply-accumulate operations with no active power \cite{Xu_Tan_Corcoran_Wu_Boes_Nguyen_Chu_Little_Hicks_Morandotti_etal._2021}. 

However, the ability to leverage favorable scaling dynamics and other properties that may be offered by novel devices rests upon the ability to meaningfully manipulate, transport, and copy the information they contain on a large scale. Without this ability, devices cannot be integrated into practical computing systems. The history of computing contains a multitude of alternative visions that were unsuccessful because of a lack of scalability, incompatibility with mainstream applications, difficulty in programming, or other issues \cite{Hooker_2021}. Compatibility with a broader model of computing provides one reassurance that a novel device can meet the end goal of being integrated within a system that provides utility to a significant user base. 

Today, systems utilizing electronic digital logic dominate mainstream computing systems for a variety of reasons: deterministic operation, interoperability, and the incredible scaling dynamics of electronic silicon technologies that persisted for decades \cite{Thompson_2017}. However, the slowdown in progress in these systems has motivated a re-exploration of alternative computing models, such as processors specialized toward machine learning \cite{thompsonDeclineComputersGeneral2018}. Many of these systems remain based on digital logic and Turing-complete designs \cite{laics}. This approach provides theoretical guarantees on their capabilities and reliability but limits the scope of devices and representations that may be leveraged within these systems. 

In contrast, analog computing can directly represent rational and real values with physical properties that vary continuously. This ability allows for intrinsic hardware trade-offs among power efficiency, storage density, and accuracy \cite{Yu_Jiang_Huang_Peng_Lu_2021, Amirsoleimani_Alibart_Yon_Xu_Pazhouhandeh_Ecoffey_Beilliard_Genov_Drouin_2020}. Analog algorithms can also offer novel solutions to NP-hard problems by leveraging real-valued dynamical systems \cite{ksat}. However, these advantages come at the cost of contending with limitations on the dynamic range of representations, innate hardware variability, probabilistic outputs, a lack of high-level programming paradigms, and other concerns. 

Oscillatory computing leverages the ability of devices that can directly represent values on an angular or ``phase'' domain at one or more frequencies \cite{Csaba_Porod_2020}. We demonstrate that by pairing this representation of information with the model of hyperdimensional (HD) computing, an analog computing model emerges that has effective high-level programming paradigms as well as clear circuit and integration primitives. We achieve this by deriving operations for an HD computing system in the complex plane that can be expressed by oscillators in the analog domain. These oscillator-based operations are implemented along with their conventional counterparts that represent phases by floating point (Figure \ref{fig:conceptual}). This dual implementation is applied to demonstrate that the two approaches can perform comparably, both in the basic operations of the HD computing system and in practical algorithms that include graph compression, factorization of composite symbols, and neural network classification driven directly via analog signals. This provides a new avenue for the design and integration of effective oscillator-based computation systems.

\section{Phase Vectors}

HD computing defines computing methods that employ phenomena encountered in high-dimensional spaces to address a variety of transformations and applications \cite{Kanerva_2009, Neubert_Schubert_Protzel_2019, Kleyko_Rachkovskij_Osipov_Rahimi_2021a, Kleyko_Rachkovskij_Osipov_Rahimi_2021b}. The fundamental units of information (symbols) in HD systems consist not of a single value but rather of a long vector of 
values---generally, at least hundreds of values. Previously, Kleyko et al. proposed that HD computing provides a rich model for information processing requiring operations that can be expressed by many emerging devices \cite{Kleyko_Davies_Frady_Kanerva_Kent_Olshausen_Osipov_Rabaey_Rachkovskij_Rahimi_etal._2021}. We build off this proposition by employing the Fourier Holographic Reduced Representation (FHRR) system, an HD system in which each symbol is defined by a vector of angular values or, equivalently, a vector of complex values on the unit circle (\ref{eq1}) \cite{Plate_1995, Plate_2003}.

\begin{equation}
\bm{\phi} \in FHRR = [\phi_1, \phi_2, \ldots, \phi_n] = [e^{i \phi_1}, e^{i \phi_2}, \ldots, e^{i \phi_n}]
\label{eq1}
\end{equation}

\begin{figure}[htbp]
    \centerline{\includegraphics[width=0.5\textwidth]{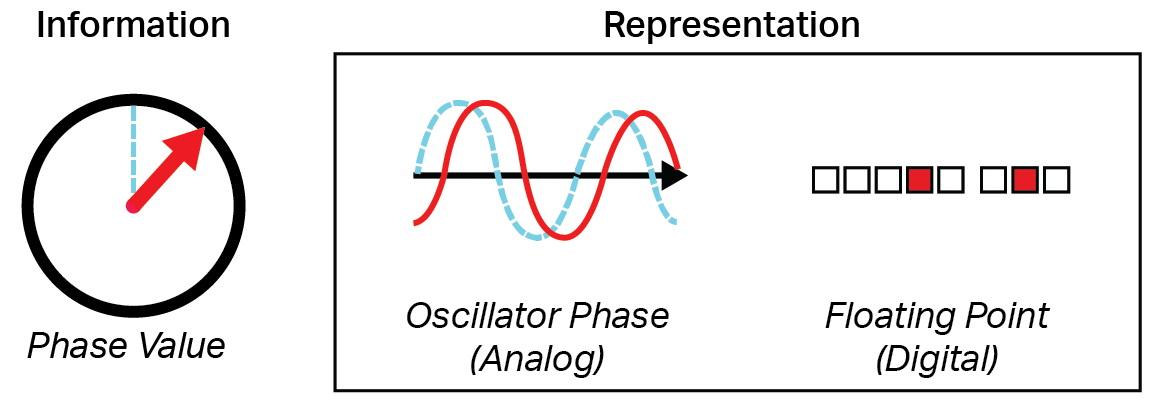}}
    \caption{Phase values can be thought of as angles subtended around a unit circle (left). This information is the basis of the hyperdimensional (HD) computing system, the Fourier Holographic Reduced Representation (FHRR), which is applied in this work. We demonstrate that this HD computing system can be effectively executed using representations of phase values expressed via oscillatory devices (middle) and digital floating-point representations (right).}
    \label{fig:conceptual}
\end{figure}

This HD system establishes a domain on which any oscillatory device or circuit can express values (Figure \ref{fig:conceptual}). Three main operations are used to manipulate and compare the information stored within an HD symbol: similarity, bundling, and binding. In this section we detail these operations, derive the methods through which oscillators can be applied to store HD symbols, and calculate these operations.

\subsection{Similarity}
As a distance can be defined between two real numbers on the ``number line,'' so can a distance be defined between two HD symbols. This distance can be found by first computing the “similarity” between two symbols \cite{Neubert_Schubert_Protzel_2019, Kleyko_Davies_Frady_Kanerva_Kent_Olshausen_Osipov_Rabaey_Rachkovskij_Rahimi_etal._2021}. Computing this similarity requires finding the phase difference between each element of the two symbols, measuring its cosine distance, and taking the average value of these distances (\ref{eq2}). Vectors of phases with $n$ elements are represented by $\bm{\phi}_0$ and $\bm{\phi}_1$. 

\begin{equation}
similarity(\bm{\phi}_0, \bm{\phi}_1) = \frac{1}{n} \sum_{k=1}^n cos(\bm{\phi}_{0,k}, \bm{\phi}_{1,k})
\label{eq2}
\end{equation}

The distance between two symbols can then be measured by subtracting the similarity value from 1---the maximum similarity (closeness) between two symbols (\ref{eq3}). 

\begin{equation}
distance(\bm{\phi}_0, \bm{\phi}_1) = 1 - similarity(\bm{\phi}_0, \bm{\phi}_1)
\label{eq3}
\end{equation}

\begin{figure}[htbp]
\centerline{\includegraphics[width=0.5\textwidth]{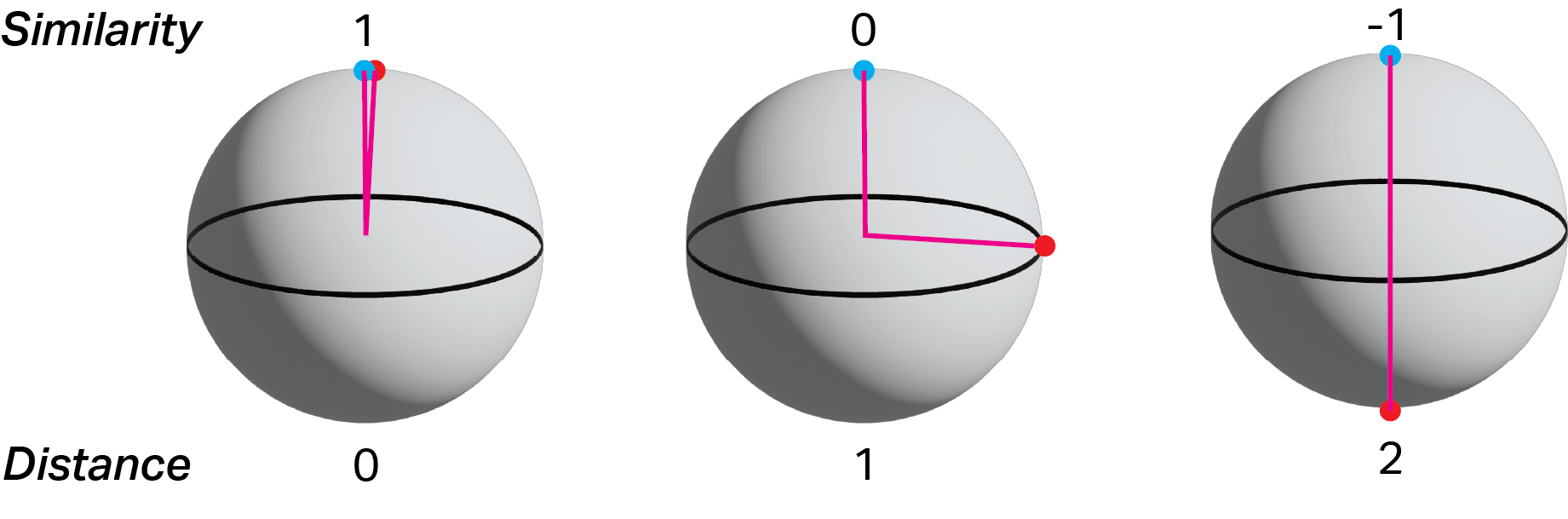}}
\caption{Illustration of similarity between HD symbols. Each FHRR symbol can be embedded on the surface of a hypersphere. Coincident points are similar (in phase), orthogonal points are dissimilar (out of phase), and points on opposite ends of the hypersphere are opposing (antiphase). Concentration of measure implies that as the number of dimensions in a symbol increases, the expectation of similarity between two random symbols tends toward zero (center diagram).}
\label{fig2}
\end{figure}

Two symbols that are identical therefore have a similarity of 1 and distance of 0. In the FHRR, this implies that each pair of angles between the two symbols is exactly in phase. Alternatively, when on average the angles are orthogonal, the similarity between the symbols will be 0, and their distance is 1. Symbols that contain angles exactly opposing one another have a similarity of -1 and distance of 2 (Figure \ref{fig2}). A fundamental phenomenon leveraged by HD computing is ``quasi-orthogonality'' of random symbols; that is, symbols constructed from randomly selected values will almost always have a similarity close to zero \cite{Kleyko_Davies_Frady_Kanerva_Kent_Olshausen_Osipov_Rabaey_Rachkovskij_Rahimi_etal._2021}.

\subsection{Bundling}

\begin{figure}[htbp]
\centerline{\includegraphics[width=0.5\textwidth]{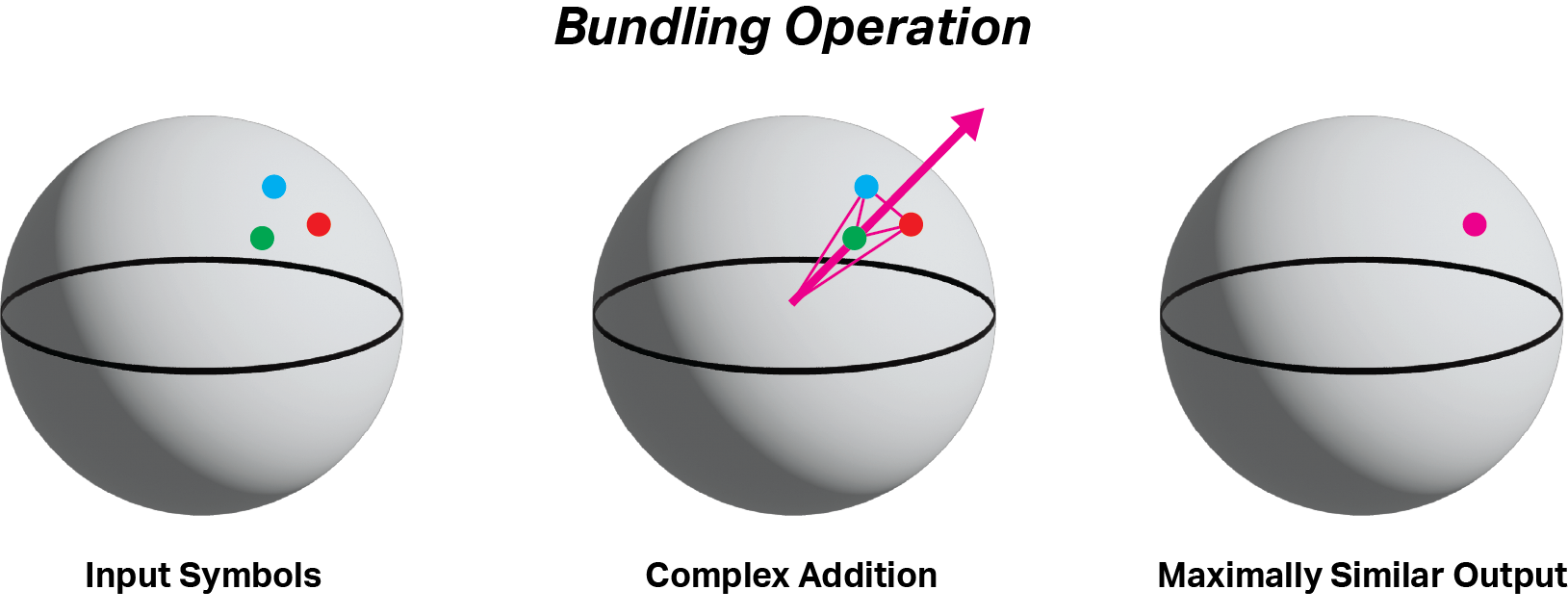}}
\caption{Illustration of the FHRR bundling operation. A set of points is combined through addition in the complex domain. The argument of this vector is taken to find a point on the surface of a hypersphere that is maximally similar to the set of inputs.}
\label{fig3}
\end{figure}

``Bundling'' is an operation in which multiple symbols are taken as an input and reduced into a single output that is maximally similar to its inputs (Figure \ref{fig3}). Creating a symbol that is maximally similar to a set of inputs requires that the cosine distance between phases in the inputs and the output phase be minimized. This is done for $m$ input symbols by converting each phase to an explicitly complex value and taking the argument of their sum (4). 

\begin{equation}
\bm{\phi}' = bundle(\bm{\phi}_0, \bm{\phi}_1, \ldots, \bm{\phi}_m) = arg(\sum_{k=1}^{m}e^{i\bm{\phi}_k})
\label{eq4}
\end{equation}

Bundling can be used as a form of lossy compression, for instance representing a set of input symbols (“apple, banana, grape”) as a single output (“fruits”). The constructed output is maximally similar to each of the input symbols. In general, symbols can be bundled until the similarity of the output to each of the inputs degrades to the level seen between random symbols. If inputs already retain a degree of similarity to one another, many may be effectively bundled together. Dissimilar inputs may also be bundled together to produce a similar output, but fewer can be included before the level of similarity degrades to the level expected between random symbols (near zero) \cite{Kleyko_Rachkovskij_Osipov_Rahimi_2021a, Neubert_Schubert_Protzel_2019}.

The inverse operation of bundling can also be carried out by scaling the output vector and removing the inputs that were added to reach it. However, this requires retrieving the other vectors that were combined to produce the output (\ref{eq5}).  

\begin{equation}
\begin{split}
\bm{\phi}_0 = unbundle(\bm{\phi}', (\bm{\phi}_1, \ldots, \bm{\phi}_m)) \\ = arg[m \cdot e^{i\bm{\phi}'} + \sum_{k=1}^{m}(-1+0i) \cdot e^{\bm{\phi}_k}]
\label{eq5}
\end{split}
\end{equation}

\subsection{Binding}

\begin{figure}[htbp]
\centerline{\includegraphics[width=0.5\textwidth]{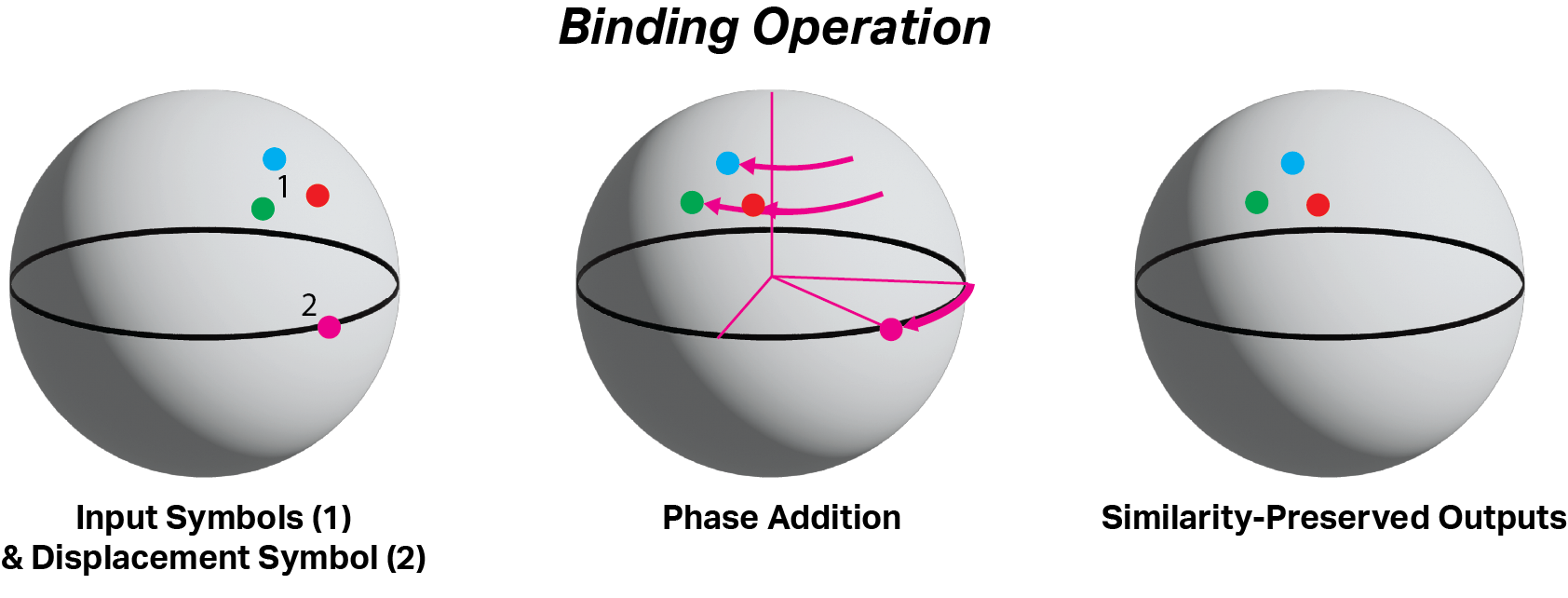}}
\caption{Illustration of the FHRR binding operation. A set of inputs (1) is displaced by a second symbol (2). The angular values of this symbol are added to the set of inputs in the phase domain, shifting each equally. This “orbits” the points around the surface of the hypersphere, maintaining each symbol’s similarity to one another but making them dissimilar to their original values.}
\label{fig4}
\end{figure}

``Binding'' is an operation in which a set of inputs can be displaced by another singular vector (Figure \ref{fig4}). Carrying out the binding operation is simple in the phase domain; one phase value is offset by another---“adding” the two angles together (\ref{eq6}). The resulting outputs can be dissimilar to the inputs, but the relationship of similarities between inputs is preserved. 

\begin{equation}
\bm{\phi}'' = bind(\bm{\phi}_0, \bm{\phi}_1) = \bm{\phi}_0 +  \bm{\phi}_1
\label{eq6}
\end{equation}

Binding can also be inverted through subtraction of phase angles (\ref{eq7}).

\begin{equation}
\bm{\phi}_0 = unbind(\bm{\phi}'', \bm{\phi}_1) = \bm{\phi}'' - \bm{\phi}_1
\label{eq7}
\end{equation}

Binding can be employed to``compose'' symbols to many ends; for instance, to represent a family member, a person’s name can be bound with the person's role, such as “Alice” with “Sister.” Multiple relatives can then be combined through bundling to produce a compressed family tree. 

\subsection{Summary}

By applying the operations of bundling and binding, a large variety of useful applications can be expressed, including finite state machines, data structures, and string processing \cite{Kleyko_Rachkovskij_Osipov_Rahimi_2021b}. By supplementing these operations with standard linear algebra, HD computing can also be used to implement neural networks \cite{Olin-Ammentorp_Bazhenov, Olin-Ammentorp_Bazhenov_2022, Bybee_Frady_Sommer_2022}. Beyond conventional computing, HD computing also offers the ability to compute multiple queries in superposition and efficiently encode and compute with sparse information \cite{Frady_Sommer_2019, Kent_Frady_Sommer_Olshausen_2019, Frady_Kent_Olshausen_Sommer_2020}. Previous work has demonstrated that neuromorphic implementations of HD computing systems are viable \cite{Orchard_Jarvis_2023}. In this work we advance this understanding by firmly establishing the connection between FHRR operations and a temporal, oscillator-based execution method.

\section{Oscillators}

Numerous natural systems exhibit oscillatory behavior: fireflies, the human heart, and pendulum clocks are all oscillators. Utilizing power from an internal source, they maintain regular, self-sustaining periodic activity that follows a regular ``limit cycle'' around an attractor in their phase space \cite{Pikovsky_Rosenblum_Kurths_2001}. In a pendulum clock, the motion of the pendulum continuously converts the oscillator’s energy from potential to kinetic and back again. The interchange between these two domains can be expressed as a system of linear differential equations (\ref{eq8}), where $x$ and $y$ represent variables such as the position and angular velocity of a pendulum, $b$ represents the ``damping'' or loss of the oscillator, and $\omega$ represents its angular frequency.

\begin{equation}
\frac{dx}{dt} = -\omega y + bx, \frac{dy}{dt} = \omega x + by
\label{eq8}
\end{equation}

These two separate variables can be transformed into a single, complex-valued argument or ``state'' that evolves through time given a single equation (\ref{eq9}), where $Z$ is the complex state and $i$ is the imaginary unit \cite{Izhikevich_2001}.

\begin{equation}
\frac{dZ}{dt} = (b + i \omega)Z
\label{eq9}
\end{equation}

Biological neurons can also be viewed as fundamentally oscillatory devices, since they exhibit frequency resonances and other complex behaviors. For this reason, Equation 9 is used as the basis for the resonate-and-fire neuron model, where the real part of the state represents the membrane current of a neuron and the imaginary part represents its membrane potential \cite{Izhikevich_2001}.

When the damping value $b$ of (\ref{eq9}) is small compared with the angular frequency $\omega$, the state of the system can be approximated over a limited period of time by Euler’s formula, (\ref{eq10}), which traces a circle around the origin of the complex plane. The value $\phi$ represents the ``starting phase'' of the oscillator that arises as an integration constant. 

\begin{equation}
Z = e^{i(\omega t + \phi)}
\label{eq10}
\end{equation}

The state of this oscillator at an instant in time can then be described as a single value---its angular position in the complex plane, known as its instantaneous phase. This instantaneous phase changes continuously with time as the oscillator’s state evolves. 

The evolution of this system in the complex plane and its ability to maintain and manipulate phase values suggest that it can provide a basis for HD computing with the FHRR system. In the following sections we prove that this intuition is correct by providing formulas to calculate HD operations directly via operations in the complex plane. 

\section{HD Operations via Oscillators}
\subsection{Representing Phase}

To utilize oscillatory elements to represent and compute with the phase values used in the FHRR HD computing system, one must encode these values in a way that remains invariant with time. To  this end, two oscillators with the same fundamental frequency can be used. Stable phase values can then be encoded, not in each oscillator’s instantaneous phase, but in their relative phase, which remains constant with time (\ref{eq11}, Appendix \ref{sec:sup_proof_1}). The time invariance of the relative phase of frequency-locked oscillators endows these systems with the ability to compute using the FHRR computing system. In other words, ideal frequency-locked oscillators can maintain the differences between their starting phases through time.

\begin{equation}
\frac{\partial}{\partial t}\angle (Z_1(\phi_1) - Z_0(\phi_0)) \propto \frac{\partial}{\partial t} \| Z_1 - Z_0 \| = 0
\label{eq11}
\end{equation}

We now derive the methods through which systems of oscillators storing relative phase values can carry out the computations necessary for the FHRR HD computing system: similarity, bundling, and binding. Augmenting previous works, we demonstrate that these calculations can be directly carried out using oscillators in the complex plane without separate integrating elements that require alternative dynamics \cite{Orchard_Jarvis_2023}. 

\subsection{Similarity}

\begin{figure}[htbp]
\centerline{\includegraphics[width=0.5\textwidth]{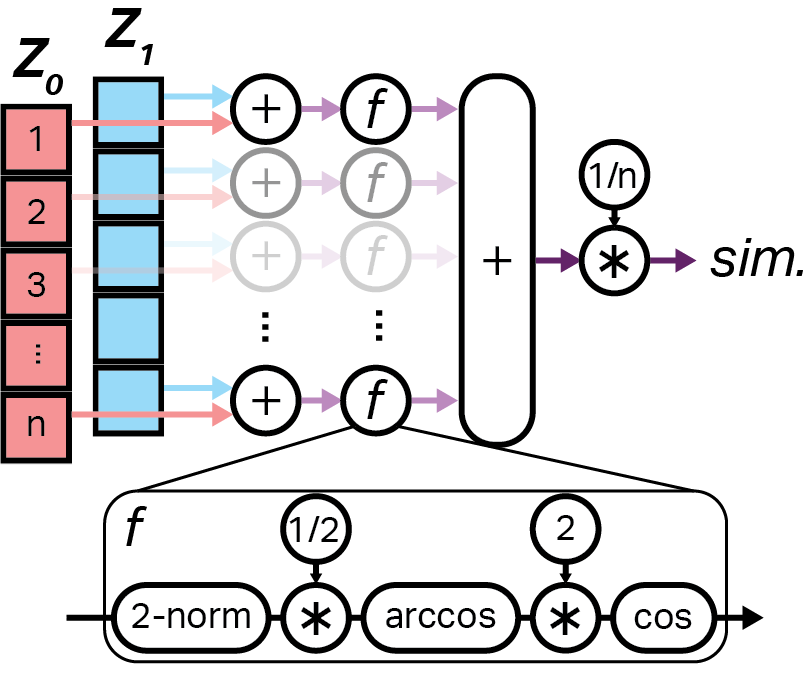}}
\caption{Illustration of the similarity operation executing on two HD symbols encoded into oscillators (blue and red squares, left). The complex values from two paired oscillators are superimposed and transformed. This operation is carried out for each pair of values between two HD symbols and reduced through averaging.}
\label{fig5}
\end{figure}

Phases can be encoded into complex state of oscillators by employing them as the starting phase of each oscillator (\ref{eq10}). As previously shown, the value of interest (the relative phase between each oscillator) will remain constant through time. Recalculating the explicit phase and computing the cosine difference between these phase angles appear necessary in order to compute similarity. However, this operation can be achieved more simply by superimposing the complex state of pairs of oscillators (Appendix \ref{sec:sup_proof_2}), where $\bm{Z}_0$ and $\bm{Z}_1$ represent vectors of complex oscillator states with $n$ elements each. This interference between the oscillators is proportional to the cosine similarity between their phase values (\ref{eq12})(Figure \ref{fig5}). 

\begin{equation}
sim.(\bm{Z}_0, \bm{Z}_1) = \frac{1}{n} \sum^{n}_{i=1} cos \left[ 2 \cdot arccos \left( \frac{\| \bm{Z}_{0,i} + \bm{Z}_{1,i} \|}{2} \right) \right]
\label{eq12}
\end{equation}

\subsection{Bundling}

\begin{figure}[htbp]
\centerline{\includegraphics[width=0.5\textwidth]{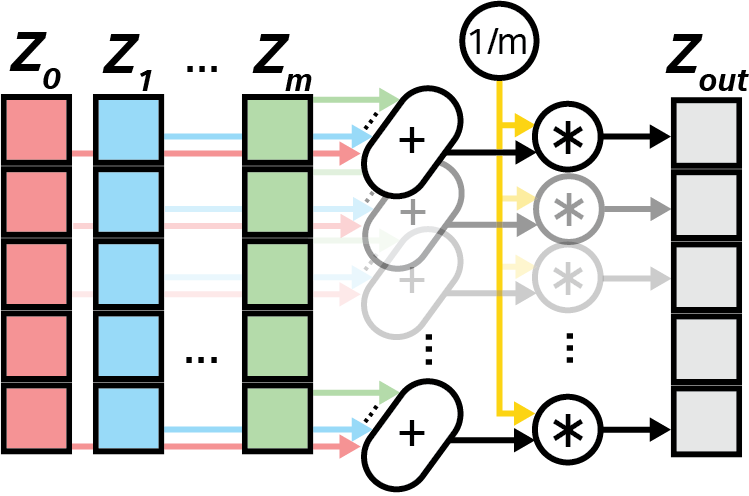}}
\caption{Illustration of the bundling operation executing on an arbitrary number of HD symbols encoded into oscillators (colored squares, left). The complex values from each oscillator in the same row are reduced through summation and normalized by the number of HD symbols to produce the output (right).}
\label{fig6}
\end{figure}

Bundling consists of a normalized superposition carried out in the complex plane. Since oscillators already encode phases in the complex domain, taking a normalized sum of their $m$ complex states carries out the bundling operation (\ref{eq13})(Figure \ref{fig6}). 

\begin{equation}
bundle(\bm{Z}_0, \bm{Z}_1, \ldots, \bm{Z}_m) = \frac{1}{m} \sum^{m}_{k=1} \bm{Z}_k
\label{eq13}
\end{equation}

\subsection{Binding}

Carrying out binding in the complex domain via the interaction of oscillators is more involved than the previous two operations, where superposition sufficed to compute the necessary values. 
In the case of binding, multiplication of complex values becomes necessary as well as the inclusion of a reference state (Appendix \ref{sec:sup_proof_3}). This reference defines what complex state currently represents the angle 0, allowing a displacement vector representing the phase of the second argument to be calculated. This allows the potential of the first argument to be ``rotated'' by the appropriate angle. The complex conjugate of the reference oscillator is included in this multiplication to negate the doubling of the output's frequency caused by multiplying two values (Appendix \ref{sec:sup_proof_3}). 

\begin{equation}
    \begin{split}
            bind(Z_0, Z_1, Z_{ref})= Z_0 + Z_0 \cdot (Z_1 - Z_{ref}) \cdot \overline{Z_{ref}}
    \end{split}
\label{eq14}
\end{equation}

\begin{figure}[htbp]
\centerline{\includegraphics[width=0.5\textwidth]{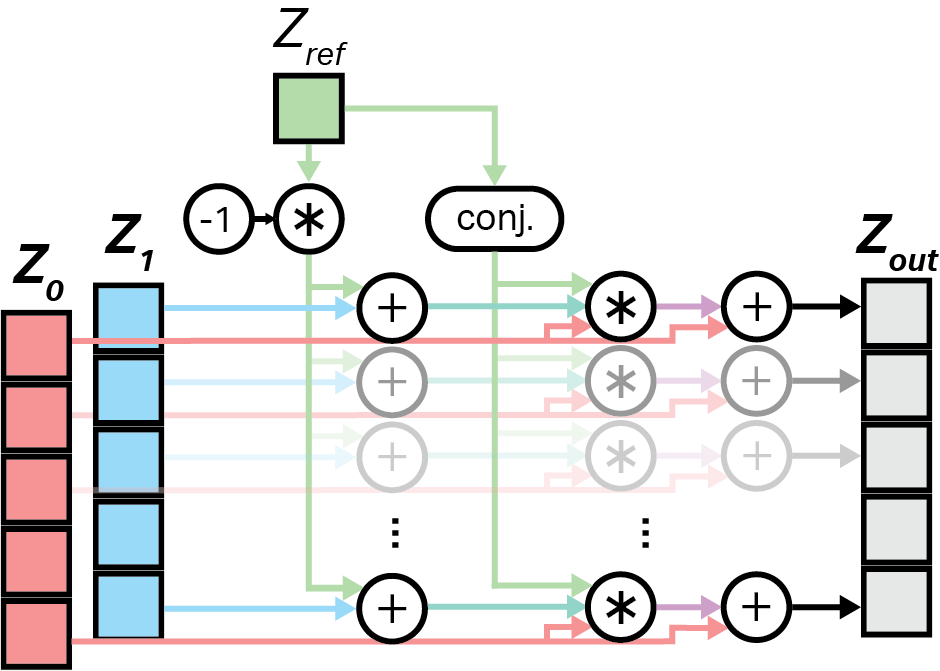}}
\caption{Illustration of the binding operation executing on two HD symbols encoded into oscillators (blue and red squares, left). In order to decode the absolute phase values stored in these oscillators, a reference oscillator representing the angle 0 is included (green square, top). The negation and complex conjugate of this reference oscillator are used to  decode the absolute phase of the second input and  maintain a frequency-matched output oscillation.}
\label{fig7}
\end{figure}

\subsection{Inverse Operations}

Since bundling takes place in the complex domain, its inverse operation can be found by negating the symbol that was added in the complex plane. Negation can be carried out by reflecting a complex value through the origin of the complex plane or adding $180^\circ$ to the equivalent phase value (\ref{eq15}).

\begin{equation}
- \bm{Z} = [-1 \cdot \bm{Z}_i] = [\bm{\phi}_i + \pi], i \in 1:n
\label{eq15}
\end{equation}

As before, the symbols involved in the original operation must be known in order to accurately invert the bundling operation.

\begin{equation}
unbundle(\bm{Z}_0 , (\bm{Z}_1, \bm{Z}_2, \ldots)) = m \cdot \bm{Z}_0 + \sum^{m}_{k=1} -1 \cdot \bm{Z}_k
\label{eq16}
\end{equation}

Conversely, binding takes place in the phase domain. To invert the binding operation, the negative of a given phase value must be added. In order to invert the phase of a complex number, it is reflected across the axis in the complex plane representing real numbers, conventionally referred to as taking the complex conjugate (\ref{eq17}). 

\begin{equation}
-\bm{\phi} = [-1 \cdot \bm{\phi}_i] = [\bar{Z_i}], i \in 1:n
\label{eq17}
\end{equation}

As previously, this defines a vector that can be used to displace the first argument ``back'' to its original value. However, complex conjugation causes the oscillator to rotate ``counterclockwise.'' The frequency correction applied by the final multiplicative term of (\ref{eq14}) must therefore be applied in the ``opposite'' direction to produce the correct output rotation.

\begin{equation}
    unbind(Z_0, Z_1, Z_{ref})= Z_0 + Z_0 \cdot \overline{(Z_1 - Z_{ref})} \cdot Z_{ref}
\label{eq18}
\end{equation}

\subsection{Transmission}

The ability to compute requires the capability to transmit values between separate parts of a 
calculation. This represents a challenge for analog computations, where transmission irreversibly distorts the values being communicated. In the case of computing via linked oscillators, the accurate transmission of many complex, time-varying analog signals makes scaling these systems to the scale of many-valued vectors required by HD computing challenging. 

One way to sidestep this issue rests on the fact that the information that must be communicated is not the entire complex state of an oscillator but just the argument of this value---its phase. This can be transmitted via impulses that are sent when an oscillator’s instantaneous phase reaches a certain value, such as 0, when the oscillator’s state is entirely real.

\begin{equation}
I(t,Z) = \delta(arg (Z(t)))
\label{eq19}
\end{equation}
\noindent
Here $\delta$ represents the Dirac delta function. These impulses may also be used to communicate an arbitrary phase from an external source.

\begin{equation}
I(t, \phi_{ext}) = \delta(arg (e^{i\omega t} - \phi_{ext}))
\label{eq20}
\end{equation}

This sparsely communicates the information necessary to represent a phase value. These pulses can serve as a source of excitement for oscillators, with real-valued current impulses causing them to resonate with the input current. 

\begin{equation}
\frac{dZ}{dt} = (b + i\omega)Z + I(t)
\label{eq21}
\end{equation}

This operation allows for systems of oscillators to communicate and synchronize via temporally sparse impulses or ``spikes,'' where the exact amplitude matters less than the timing. 

\section{Demonstrations}

We now demonstrate HD computing operations carried out via two separate methodologies: (1) a standard implementation in which floating points directly represent phase angles and (2) a system in which phase angles are encoded into binary pulses (“spikes”) that excite oscillators in the complex domain. To simulate these temporal systems, we solve (\ref{eq21}) numerically through time using the constants $b=-0.2$ and $\omega = 2\pi$ Hz.

The performance of each operation executing via the interaction of oscillators is tested by computing its output between all possible pairs of values on the domain of angular values ($[0,2\pi]$) and comparing it with the standard, atemporal floating-point implementation. All demonstrations are made by using our Julia software package PhasorNetworks.jl, which implements both execution techniques.  

\subsection{Similarity}

\begin{figure}[htbp]
\centerline{\includegraphics[width=0.5\textwidth]{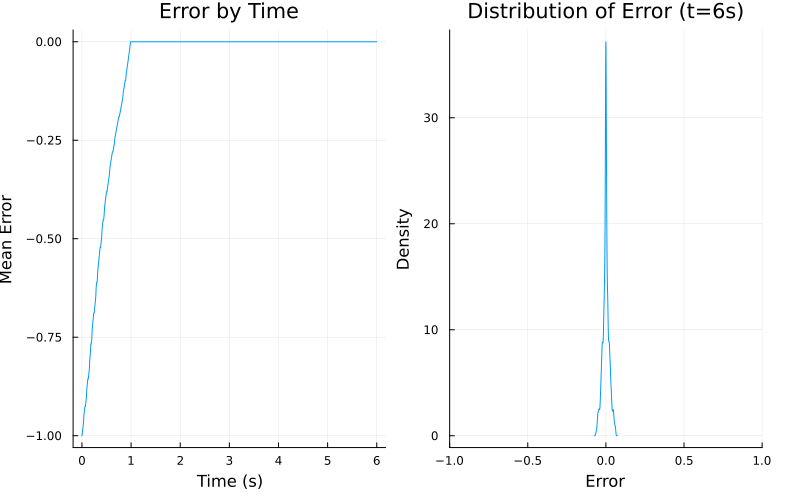}}
\caption{The mean amount of error between values calculated by a standard similarity operation and one operating via the interference of oscillators by time is shown (left). The distribution of these error values at the end of the simulation is narrow and centered on zero (right).}
\label{fig8}
\end{figure}

The difference in values output by the standard similarity function (\ref{eq2}) and the oscillator-based similarity function (\ref{eq12}) through time is displayed in Figure \ref{fig8}. These similarity values fall in the real-valued domain $(-1,1]$. The error in similarity values decreases rapidly through time as input pulses transmitted to oscillators cause them to resonate with the correct phase. The complex potentials of these resonating oscillators then interfere to obtain the desired similarity values, with differences from the floating-point method remaining small and distributed around zero. 

\subsection{Bundling}

\begin{figure}[htbp]
\centerline{\includegraphics[width=0.5\textwidth]{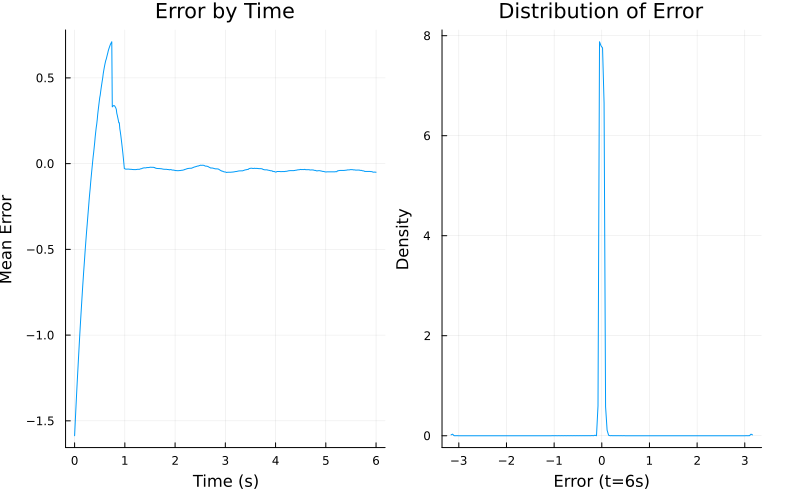}}
\caption{The mean amount of error between phase values produced by a standard bundling operation and one operating via the superposition of oscillator states is shown. This error decreases through time as the oscillators resonate with their input signals (left). The distribution of this error at the final simulation step is narrow and centered on zero (right).}
\label{fig9}
\end{figure}

The difference in phase values output by the standard bundling function (\ref{eq4}) and the oscillator-based bundling function (\ref{eq13}) is displayed in Figure \ref{fig9}. Again, as the oscillators begin to resonate with their driving input signals, the superposition of their states produces the necessary outputs with small errors.

\subsection{Binding}

\begin{figure}[htbp]
\centerline{\includegraphics[width=0.5\textwidth]{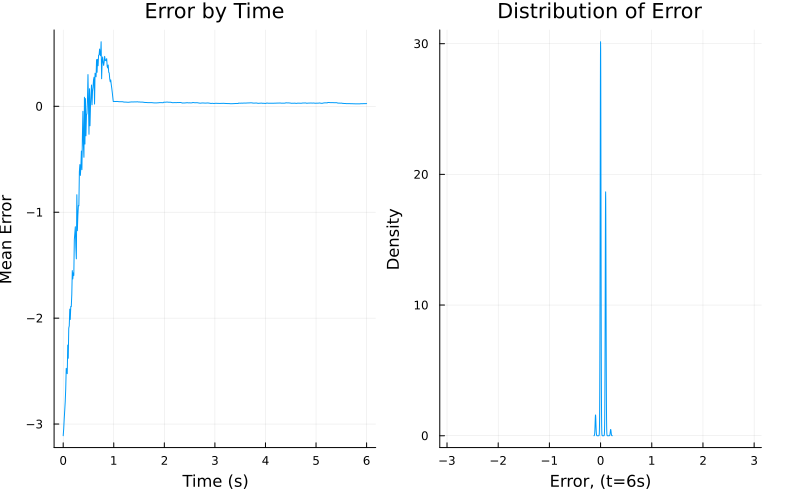}}
\caption{The mean amount of error between phase values produced by a standard binding operation and one operating via the multiplication of oscillator states is shown. As in the previous operations, this error decreases through time (left) and becomes narrowly centered on zero (right).}
\label{fig10}
\end{figure}

We also carried out this comparison  for the binding operation, with the standard method (\ref{eq5}) being used as the baseline against which the oscillator-based method (\ref{eq14}) is compared. The error of the oscillator-based method again decreases through time and becomes small and centered on zero (Figure \ref{fig10}).

\section{Algorithms}

The previous experiments characterized the efficacy of FHRR computations implemented via oscillators on individual pairs of angles---these are the individual numerical operations constituting the FHRR computing system. Next, we demonstrate the efficacy of these individual operations as they scale to full hyperdimensional vectors that are manipulated to carry out practical, multistage algorithms that can execute effectively on analog systems.

\subsection{Graph Compression}
\subsubsection{Overview}
The first algorithm demonstrates the compression and reconstruction of graphs via HD symbols, accomplished by representing all edges contained within the graph with a single symbol \cite{Poduval_Alimohamadi_Zakeri_Imani_Najafi_Givargis_Imani_2022, Gayler_Levy}. Completing this task requires composing all previously demonstrated functions: binding, bundling, and similarity (as well as their inverse operations). 

A series of undirected graphs with no self-loops were constructed via the Erdős-–Rényi model. These graphs contain $n$ nodes and a variable number of edges selected out of all possible pairs with probability $p$. Increasing $p$ thus leads to an increase in the expected number of edges per graph. 

These graphs can be represented in an HD space by selecting a set of random HD symbols, where each symbol uniquely corresponds to one node (\ref{eq22}). In this experiment each graph contained 25 nodes, and each symbol --- corresponding to a node --- contained 1,024 phase values. 

\begin{equation}
node_i := \bm{\phi}_i
\label{eq22}
\end{equation}

In order to describe each edge in the graph, the symbols of the nodes adjacent to the edge can be bound to produce a new symbol (\ref{eq23}).

\begin{equation}
edge_i = [node_j \Leftrightarrow node_k] := bind(\bm{\phi}_j, \bm{\phi}_k)
\label{eq23}
\end{equation}

This creates two sets of symbols: one representing the graph’s nodes and the other representing its edges. The HD properties of these symbols makes it highly probable that each node symbol will be dissimilar to all others; and, being derived from them, all edge symbols will be dissimilar to both the node symbols and other edge symbols. The amount of space needed to store the set of edges can be reduced via bundling to a single symbol (\ref{eq24})(Listing 1). 

\begin{equation}
\bm{G} := bundle(edge_1, edge_2, \ldots, edge_m)
\label{eq24}
\end{equation}

\begin{lstlisting}[caption = Compressing a graph's edge structure into a single HD symbol]
function compress_edges(graph):
#represent each node with an HD symbol
HD_nodes = rand_syms(length(graph.nodes))
HD_edges = []

for edge in graph.edges:
  #produce a unique representation 
  # of each edge
  source_index, destination_index = edge
  HD_edge = bind(HD_nodes[source_index], 
    HD_nodes[destination_index])

  append!(HD_edges, HD_edge)
end

#reduce the list of edges 
# to a single symbol via bundling
HD_graph = bundle(HD_edges...)

return HD_nodes, HD_graph
\end{lstlisting}

The symbol $\bm{G}$ can be thought of as representing the list of edges in a graph in compressed form via HD operations. To reconstruct the graph’s adjacency matrix from this symbol, we use the set of symbols representing the nodes N. For each node, its corresponding symbol is unbound from $\bm{G}$, and the similarity of this product to all other nodes is calculated. This process produces a matrix of similarity values (\ref{eq25})(Listing 2). 

\begin{equation}
\begin{split}
\bm{A} = [sim.(node_j, unbind(\bm{A}, node_i))] \\ for i \in 1:n, j \in 1:n, i \neq j
\end{split}
\label{eq25}
\end{equation}

\begin{lstlisting}[caption = Reconstructing a graph's edges from a node symbols and a compressed edge symbol]
function predict_edges(HD_nodes, HD_graph):
 n_nodes = length(HD_nodes)
 adjacency = zeros(n_nodes, n_nodes)

 #unbind each node from the graph
 # symbol to produce a key
 for i in 1:n_nodes
  source_node = HD_nodes[i]
  key = unbind(HD_graph, source_node)

  #measure the similarity between this key
  # and another node
  for j in 1:n_nodes
   query = HD_nodes[j]
   adjacency[i,j] = similarity(key, query)
  end
 end

 return adjacency
end
\end{lstlisting}

\subsubsection{Results}
All values in $\bm{A}$ above a certain threshold can be accepted as predicted true edges, and the rest as predicted false edges. This  classification  can be compared against the ground truth and has an overall performance measured by the area under a receiver operating curve (AUROC). 

\begin{figure}[htbp]
\centerline{\includegraphics[width=0.5\textwidth]{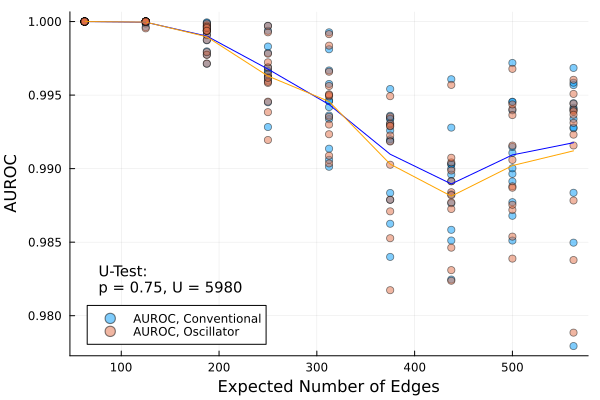}}
\caption{The performance of a classifier used to reconstruct graphs from a compressed representation was measured against an increase in the probability of edges in an Erdős--Rényi random graph. As the encoding process increases in complexity, the quality of the reconstruction determined from the HD embedding suffers. Executing this operation via oscillator-based computation as opposed to a floating-point implementation showed more significant decreases more quickly, demonstrating that small changes in the accuracy of HD operations can have significant impacts on performance. }
\label{fig11}
\end{figure}

This compression task increases with difficulty as more edges are included in the graph until a point of maximum entropy is reached. This is reflected in the trend of the performance of the reconstruction with increasing probability of edges in the graph for the conventional floating-point implementation (Figure \ref{fig11}). The oscillator-based implementation produces values that produce orderings indistinguishable from the floating-point implementation (U-test, $p = 0.75, U = 5980$) and display an identical trend (Figure \ref{fig11})

\subsection{Symbolic Factorization}
\subsubsection{Overview}
One advantage of HD computing is the ability to compose new symbols from an extant vocabulary through binding. For instance, in the previous example, symbols representing nodes were manipulated via binding to represent edges between nodes. However, to ``decipher'' these new representations, one may have to ``factor'' a symbol into its original vocabulary. This factoring can be accomplished by brute force, as in the previous example, where all potential combinations are examined one by one. As the number of components bound to create a symbol increases, however, the factoring problem becomes exponentially more complex.

Resonator networks were proposed as a solution to this challenge. These networks leverage the ability of HD computing systems to carry out queries ``in superposition'' \cite{resonator1,resonator2}. Essentially, by bundling an entire vocabulary together into a ``guess'' at a factor, all combinations of factors are tested simultaneously. Improving this guess over several iterative steps creates a dynamical system that is likely to converge to the correct factorization. 

In order to formally define this problem, a symbol can be composed by binding elements selected from different sets of symbols. For instance, one set of symbols $\textbf{A}$ may define $n$ colors, and another set of symbols $\textbf{B}$ may be used to define $m$ animals. Binding can thus be used to combine two individual symbols $i$ and $j$ from each set to produce a new symbol.

\begin{equation}
    z_{i,j} = bind(\textbf{A}_{i}, \textbf{B}_{j})
\label{eq26}
\end{equation}

If the color in the set $\textbf{A}$ at index $i$ represents ``black'' and the animal at index $j$ in set $\textbf{B}$ represents ``cat,'' the meaning of the symbol $z$ is a ``black cat.'' Binding these vocabularies can represent a total of $n \cdot m$ colored animals. However, the binding operation produces an output that is dissimilar to all input symbols in sets $\textbf{A}$ and $\textbf{B}$. To decode this meaning of a symbol whose constituents are not known, one must solve the inverse problem.

\begin{equation}
    \forall z \in bind(\textbf{A}_{i \in n}, \textbf{B}_{j \in m}), i=?, j=?
\label{eq27}
\end{equation}

A resonator network makes an initial guess for each component $a'$, $b'$, created by bundling all elements in sets $\textbf{A}$ and $\textbf{B}$ to place each guess on the span of these vocabularies. This allows all potential $m \cdot n$ combinations of elements to be tested ``in superposition.'' 

\begin{equation}
    \begin{split}
        a' = bundle(a \in A...) \\ 
        b' = bundle(b \in B...) 
    \end{split}
\label{eq28}
\end{equation}

These guesses are then bound to create a prediction $z'$ that will be iterated to reconstruct $z$.

\begin{equation}
    z' = bind(a', b')
\label{eq29}
\end{equation}

Access to each of the ``guesses'' allows the contributions from other factors to be refined. In order to assess the quality of the current guesses, they are applied to mutually reconstruct other factors. For instance, in natural numbers, the composite product of 15 and prime factor of 3 allows one to infer 5 as the separate prime factor. Correspondingly, unbinding the current guess $b'$ from the reconstruction $z'$ produces a new guess, $a''$, for the true color factor $a$ independent of the current guess $a'$. The similarity of $a''$ to the elements of the set $\textbf{A}$ can then be used within bundling to create a ``cleaned-up'' factor, with the most similar elements carrying the highest weight. This independent prediction produces an improved guess, $a'''$. 

\begin{equation}
    \begin{split}
    a'' = unbind(z', b') \\
    a''' = similarity(a'', \textbf{A}) \cdot \textbf{A} \\
    \end{split}
\label{eq30}
\end{equation}

This same process is carried out for other factors---in this case $b'''$, but scaling to any number of factors. Once all factors have been refined, a new guess for the composite factor, $z''$, is produced. 

\begin{equation}
    \begin{split}
    z'' = bind(a''', b''')
    \end{split}
\label{eq31}
\end{equation}

This iterative algorithm is likely to converge quickly to the correct answer for each factor because of the pseudo-orthogonality of each factor, and is expressed via pseudocode in Listing 3 \cite{resonator1}. The above example uses two factors for clarity of notation, but it can be expanded for an arbitrary number of factors: 

\begin{lstlisting}[caption = Factoring a composite symbolic representation]
function refine_guess(guess, factor_set)
  weights = similarity(guess, factor_set)
  refined_symbol =
    bundle((weights * factor_set)...)
  return refined_symbol
end

function predict_composite(guesses...)
 composite = reduce(bind, guesses...)
end

function likely_factors(guess, factor_set)
 return argmax(
          similarity(guess, factor_set))
end
 
function factor(composite, n_iterations, 
                        factor_sets...):
  guesses = [bundle(factors...) 
             for factors in factor_sets]
  n_factors = length(factor_sets)

  for i in 1:n_iterations
    reconstruction = 
        predict_composite(guesses)
    new_guesses = 
      [refine_guess(g, 
      factor_sets[setdiff(1:n_factors, j))] 
      for j in 1:n_factors]
    new_reconstruction = 
        predict_composite(new_guesses)
    guesses = new_guesses
    reconstruction = new_reconstruction
  end

  factors = [likely_factors(guesses[i],
    factor_sets[i])
    for i in 1:n_factors]
  return factors
end
\end{lstlisting}

\subsubsection{Results}

\begin{figure}[htbp]
    \centering
    \includegraphics[width=1\linewidth]{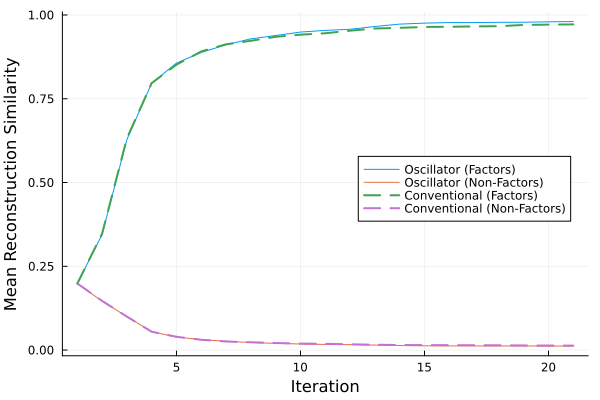}
    \caption{Average similarity of the composite symbol reconstructed by a resonator network in comparison with the factors and non-factors (n=256). Over 20 iterations of the network, solutions tend to converge quickly to a high-quality reconstruction that allows the correct factors to be identified.}
\label{fig:resonator_network}
\end{figure}

This program was implemented for HD symbols constructed from three factors, each drawn randomly from a set of 20 symbols, with each symbol containing 1,024 phase values. The resonator network was run for 20 iterations, and the predicted factors were checked for correctness against the ground truth. Within these iterations, 97\% of the floating-point-based solutions were correct, and 98\% of the oscillatory solutions were correct---a difference that is not statistically significant  (U-test, $p = 0.18, U = 298e3$). On average, the similarity of reconstructed composite symbols over time converges quickly to the correct solution for both implementations (Figure \ref{fig:resonator_network}). This situation demonstrates both the performance of resonator networks in the FHRR HD computing framework and the feasibility of operating a second algorithm via oscillatory operations. 

\subsection{Analog Neural Networks}
\subsubsection{Overview}

Neural networks have grown into one of the most important computing applications in recent years \cite{diffusion_of_ai}. The growth of energy usage for this domain alone has led to rapid growth in the energy usage and overall environmental impact of computing \cite{economist_data,environmental_impact}. As a result, technologies that can reduce the cost and impact of neural networks and associated artificial intelligence  technologies are an important topic of research. HDOC offers one pathway through which novel hardware technologies can be harnessed to address these concerns by providing an oscillatory implementation of neural networks \cite{Olin-Ammentorp_Bazhenov_2022}.

A neural network can be viewed as a successively weighting inputs and applying an activation function to produce a new output. This behavior can be accomplished in HDOC by utilizing bundling and complex multiplicationß. Given an $n$-dimensional, complex-valued input symbol $X$, a real-valued $m \times n$ matrix $W$ can be used to transform $X$ into $n$ symbols with dimensionality $m$. These $n$ symbols are bundled together to create a single, $m$-dimensional symbol $Y$.

\begin{equation}
    \begin{split}
    Y = NN_{HD}(X, W) =
    bundle(\textbf{W}_{:,1} \times X, \\
     \textbf{W}_{:,2} \times X, \\
     ... \textbf{W}_{:,n} \times X) 
    \end{split}
    \label{eq33}
\end{equation}

A multiple-layer network can be created by composing (\ref{eq33}) recursively with different values of $n$ and $m$ to achieve different ``layer'' sizes and the familiar multilayer perceptron (MLP) architecture.

\begin{equation}
    \begin{split}
        Y = MLP_{HD}(X, W_{1}, W_{2}) \\
        = NN_{HD}(NN_{HD}(X, W_{1}), W_{2})
    \end{split}
    \label{eq34}
\end{equation}

Multilayer networks can thus be constructed. Adjusting the weights $W_{i}$ to create the desired transformation can be accomplished by backpropagation through the locally continuous function (\ref{eq33}), either by analyzing the time-invariant phases or by applying adjoint gradient analysis methods through the oscillator-based solution to (\ref{eq33}) through time \cite{diff_p}. Additional architectures such as residual networks and attentional architectures may be implemented through HDOC operations such as binding and similarity \cite{Olin-Ammentorp_Bazhenov}.

For tasks such as classification, explicit phase symbols can be constructed as targets for the MLP's transformation to be trained toward. As in previous works, we use ``quadrature'' symbols as targets for $Y$, in which all elements of a symbol correspond to a class, yielding a dimensionality of $n_{c}$ equal to the number of classes. All phases in this symbol are zero except the position that corresponds to the true class $c$, which is set to a phase of $\pi / 2$.

\begin{equation}
    Y_{quad.} = \frac{\pi}{2} \cdot onehot(c, n_{c})
    \label{eq35}
\end{equation}

The distance between the output $Y$ and the correct quadrature label provides the metric for gradient descent to reduce loss to a minimum by adjusting the parameters $W_{i}$, training the HD neural network \cite{Olin-Ammentorp_Bazhenov_2022}.

\begin{equation}
    \begin{split}
    (W_{1}, W_{2}) = argmin(1 - \\
    similarity(NN_HD(X_{train}, W_{1}, W_{2})),\\
    Y_{quad}(Y_{train}, n_{c}))
    \end{split}
\label{eq36}
\end{equation}

Besides offering an alternative execution strategy for neural networks, HD implementations offer a unique advantage: an ability to bypass one or more costly analog-to-digital steps that may be required in order to calculate input or activation values for a neural network. Referring back to Equation \ref{eq21}, we observe that oscillators may be excited by an externally determined input current, $I(t)$: these inputs may be provided by an arbitrary source. We investigate the cases in which these inputs are provided ``directly,'' by applying a current source to drive the network, and ``indirectly,'' in which the currents are converted into values that drive the network externally. 

\subsubsection{The ``Smart Pixels'' Task}

\begin{figure}
    \centering
    \includegraphics[width=1\linewidth]{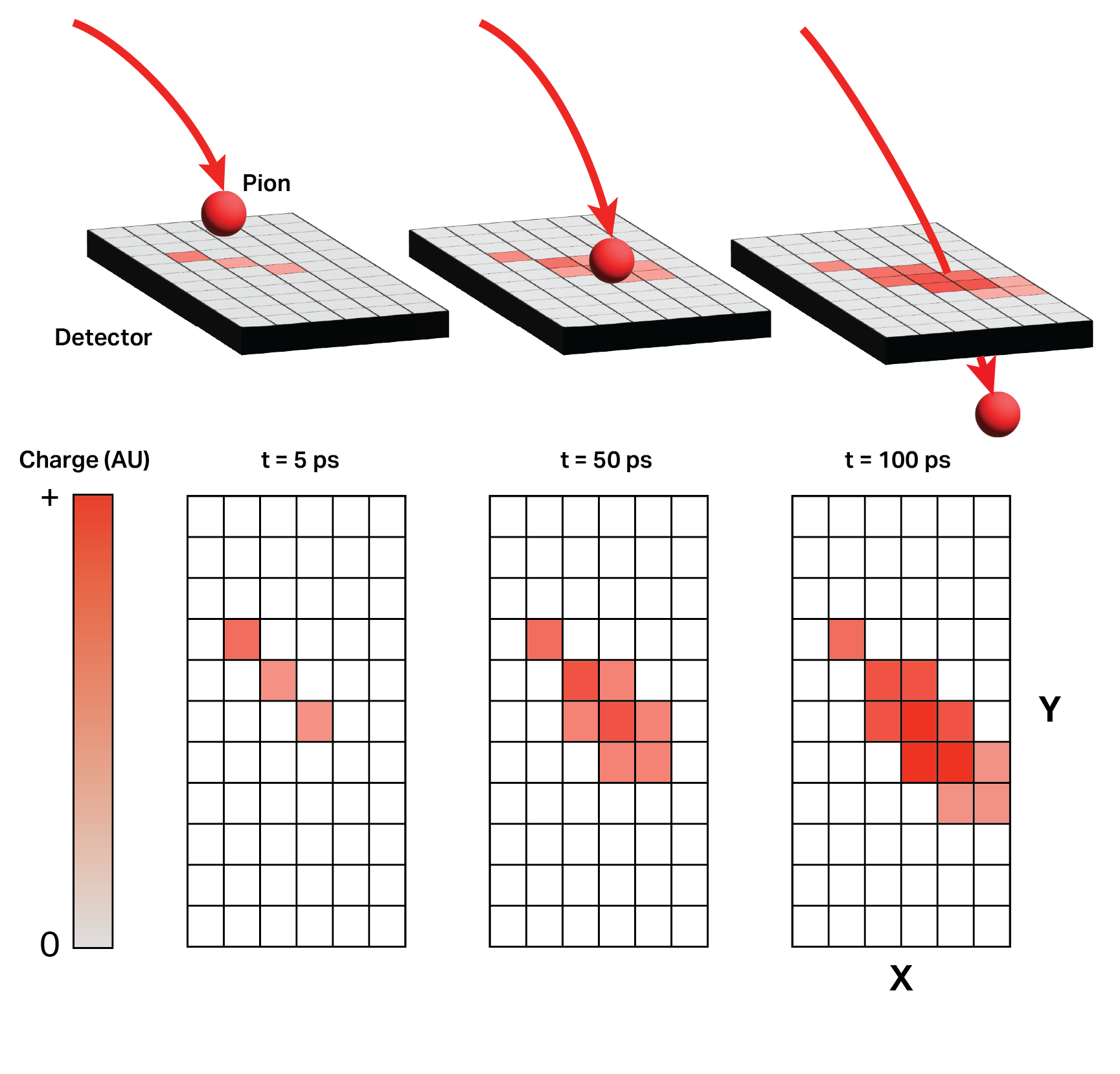}
    \caption{Illustration of the ``Smart Pixels'' dataset. Pions that transit through a thin silicon detector will cause charges to accumulate across a reverse-biased junction. The amount of charge in each pixel (red rectangles) generated varies by time as the pion (sphere) transits through the detector. Utilizing the charge generated through time and the displacement of the detector from the particle's generation point allows the pion's transverse momentum and charge to be inferred.}
\label{fig:smart_pixels}
\end{figure}

The task used for this demonstration is determining the charge of an elementary particle and whether its transverse momentum is above a threshold of 200 MeV. This is accomplished given the location and the currents generated as each particle transits through a 2-D array of silicon pixels (Figure \ref{fig:smart_pixels}) \cite{smart_pixels}. The ability to quickly determine this classification enables ``smart'' detectors that discard unneeded data from low-momentum particles.

In the indirect network the accumulated charge in each row at the end of the sampling period is scaled and converted into real or phase values that drive a standard, two-layer MLP using the ReLU activation function or an HD MLP. This information is combined with the location at which the sensor is positioned to ``catch'' the particle's transit. As with previous applications, the HD MLP is simulated by  using both conventional and oscillator-based execution.

In contrast, the direct network is an HD MLP with an initial layer of oscillators individually driven by currents produced by each row of pixels. The same location information provided to the indirect network is included in the direct inputs through a temporally-coded current pulse. As each neuron in this initial layer is excited by these input currents, it begins to resonate, producing phase values that drive the downstream two-layer HD MLP. Again, this MLP is simulated by using conventional and oscillator-based representations. The details of each transformation of charge and current into inputs for the direct and indirect neural networks are given in Appendix \ref{sec:nn_transforms}. 

\subsubsection{Results}

All networks utilize the same trainable architecture of 14 input neurons, 128 hidden neurons, and 3 output neurons. Directly driven networks contain an additional input layer of 14 oscillatory neurons; however, the weights for this layer are not trainable, and this layer simply converts input currents directly into phases. All networks were trained for 10 epochs on approximately 100,000 examples and validated against a separate set of 50,000 test examples. Since many more particles with low momentum than high momentum are present, area under the receiver operating curve (AUROC) is used to measure networks' performance rather than simple accuracy. 

\begin{table}[htbp]
    \begin{tabular}{c|c|c|c}
    Type & Type & Execution & \makecell{AUROC (max, $\mu, \sigma$, n=24)} \\ \hline
    ReLU MLP & Indirect & Conv. &  0.834, 0.829, 3.90e-3 \\
    HD MLP & Indirect & Conv. & 0.830, 0.827, 1.72e-3 \\
    HD MLP & Indirect & Osc. & 0.829, 0.825, 1.78e-3 \\
    HD MLP & Direct & Conv. & \textbf{0.862 0.855, 2.79e-3} \\
    HD MLP & Direct & Osc. & 0.857, 0.804, 3.17e-2 \\
    \end{tabular}
\caption{Summary of the experimental neural networks trained to classify pixel momentum and charge.}
\label{tab:nn_results}
\end{table}

\begin{figure}
    \centering
    \includegraphics[width=1\linewidth]{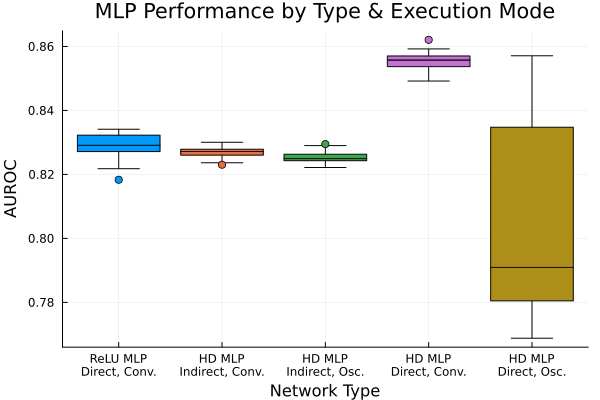}
    \caption{Comparison of the performance reached by a standard neural network vs HD neural networks under different operating conditions. Among the indirect networks, the standard ReLU-based MLP reaches the highest performance of indirectly driven networks. However, the directly driven phasor network reaches the highest performance of all compared networks, potentially due to its ability to integrate temporal information into its transformation. Unfortunately, this high performance does not always translate into the oscillatory execution mode for this network.}
\label{fig:nn_results}
\end{figure}

When the MLP networks are indirectly driven---that is, only the charge from the end of the sensor's sampling period is used as an input---performance between the networks is relatively comparable. The traditional ReLU MLP reaches higher performances than the HD MLP in either execution mode (conventional or oscillatory, Table \ref{tab:nn_results}). Nonetheless, this gain is slight, and the ReLU MLPs also produce the worst-performing indirectly driven networks (Figure \ref{fig:nn_results}); therefore, this may be a result of different initialization methods and numerics rather than limits in HD network capabilities. 

However, differences between the indirectly driven networks are small in comparison with the headway the directly driven neural networks gain in performance, which is significantly greater than that of the ReLU MLP networks (U-test, $p<1e-13, U=0$). This is likely due to the ability of the network to leverage the full temporality of the information presented to it by integrating currents through time to produce the inputs to the network. These performance gains are reliable in the conventionally executed HD MLP, but this performance is not consistently seen when executing in the oscillatory mode. This sensitivity to execution mode is not observed with the indirect HD MLP networks (Figure \ref{fig:nn_results}); and this variability in performance was sensitive to oscillatory parameters, weights, and other numerics in the underlying simulation. This suggests that the performance of these networks depends on underlying numeric conditions and could be improved by further analysis and investigation, which  we leave to future work. 

\section{Discussion}

In the preceding sections we have demonstrated that HD computing operations can be carried out via a system of linked, ideal oscillators. While this demonstration provides a proof of concept for integrating devices into useful computational circuits, the practical feasibility and relative merits of implementing these systems in hardware must be evaluated.

On the level of individual devices, all oscillators contain non-idealities that affect their performance. For instance, no two physical oscillators are perfectly matched in frequency, although they can be driven to resonate in synchrony \cite{Pikovsky_Rosenblum_Kurths_2001}. Physical oscillators also reach limits in their phase domain, such as the displacement achievable by the cantilever of a MEMS oscillator or the charge that can be stored in a resonant electronic circuit. All calculations carried out by the oscillator must be normalized within this phase space to avoid distortion or failure of components. Different components will also lead to changes in ease of fabrication, scalability, power dissipation, and accuracy of individual oscillator phases. 

A circuit-level implementation of oscillatory HD computing would require thousands of oscillators in a circuit, with series of many oscillators dedicated to represent an HD symbol. The complex-valued information in these circuits must then be accessible to carry out addition, multiplication, inversion, and conjugation: the operations required for forwards and inverse bundling and binding, as well as similarity. 

The ease of transporting this information and carrying out these operations changes greatly with the underlying physical implementation. For instance, coherent beams of light can be used to represent a phase, transported via waveguides and summed via superposition. However, ``multiplying'' two beams of light is not as simple to achieve, because of the lack of interaction between photons. While capacitive-inductive oscillators offer a similar appeal due to their ability to represent a system via two physical domains, the difficulty of integrating inductors into integrated circuits remains a challenge to their adoption in large-scale systems. Familiar electronic systems, such as CMOS, offer the highest scalability, ability to transport information, and efficient subthreshold operation. However, the use of the electronic domain alone requires at least two transistors to represent a single oscillator and explicit addition and multiplication circuits that may be achievable by simpler, physical operations in other devices. 

At the system level, the accuracy achieved by each operation must be high; our experiments suggest that small differences in the precision of HD operations can lead to significant impacts on the performance of an HD algorithm. Performance will also change with the number of values composing an HD symbol and may be applicable as a circuit-level trade-off to achieve higher accuracy using more components. Additionally, we have considered the case in which each oscillator produces a continuously valued phase; restricting these values to a finite subset would discretize the system, allowing it to carry out deterministic operations. Reducing the number of phase values would impact the overall performance of the system, but we note that one HD computing system computes digitally with only two phase values \cite{Schlegel_Neubert_Protzel_2020}.  

Furthermore, success in computing requires an operation to be easily integrable with the dominant digital ecosystem. For many emergent computing operations such as computation-in-memory, this requires a conversion between the analog and digital domains that can ultimately be costly \cite{Amirsoleimani_Alibart_Yon_Xu_Pazhouhandeh_Ecoffey_Beilliard_Genov_Drouin_2020}. Analog, oscillatory HD computing systems would require efficient, scalable means of converting phases from their analog representations into digital values that can be stored and transported to other sections of a digital computer.  

Utilizing detailed hardware models of oscillatory devices within the simulation framework is necessary in order to investigate the potential of oscillator-based HD computing systems at the device, circuit, and system level. However, the abundance of oscillatory devices such as subthreshold transistors, LC/RC circuits, MEMS cantilevers, spin torque oscillators, and coherent photonic systems suggests that there is already a rich field of devices and data that could be applied to this objective \cite{Csaba_Porod_2020}.

For many extant tasks well suited to digital logic, such as arithmetic and precise, repeatable calculations, HD computing systems are unlikely to provide competitive performance. However, for many other tasks such as interfacing with analog inputs, compressing noisy information in a structured manner, and manipulating and searching high-dimensional embeddings, HD computing is well suited to provide useful approaches \cite{Kleyko_Rachkovskij_Osipov_Rahimi_2021b, Frady_Kent_Olshausen_Sommer_2020, Imani_Rahimi_Kong_Rosing_Rabaey_2017}.

\section{Conclusion}

Developing alternative methods of computation is crucial to addressing increasing demands as development of traditional approaches continues to increase in cost. In this work we propose that hyperdimensional (HD) computing that represents information via vectors of phase angles provides a rich computational system with fundamental links to novel hardware devices. We provide proofs for the basis of this link, and we quantify via simulation the error in outputs of fundamental HD operations operating via systems of linked oscillators. We appllied these fundamental operations  to demonstrate effective graph compression, symbolic factorization, and a neural network directly excited by analog inputs. We suggest follow-up work investigating improving the accuracy of these oscillatory computations and establishing deeper links to specific classes of hardware devices. 

\section*{Acknowledgments}
We thank Andrew A. Chien, Xingfu Wu, and Angel Yanguas-Gil for discussing and revising this work and Gail Pieper for editing and proofreading. 

This work was supported by DOE ASCR and BES Microelectronics Threadwork. This material is based upon work supported by the U.S. Department of Energy, Office of Science, under contract number DE-AC02-06CH11357.

\section*{Data Availability Statement}

The data and code that support the findings of this study are openly available on GitHub, at \url{https://github.com/wilkieolin/phasor_julia/tree/main}.

\begin{appendices}
\section{Proof 1}
\label{sec:sup_proof_1}

The phase between two oscillators operating with the same frequency remains constant through time. 

\begin{equation}
\begin{split}
    \frac{\partial}{\partial t} arg(Z_1 - Z_0) \propto \frac{\partial}{\partial t} \| Z_1 - Z_0 \| \\
    \therefore \frac{\partial}{\partial t} \| e^{i(\omega t + \phi_1)} - e^{i(\omega t + \phi_0)} \| \\
    = \frac{\partial}{\partial t} \| e^{i \omega t} e^{i(\phi_1 - \phi_0)} \| \\
    = \frac{\partial}{\partial t} \| 1 \cdot e^{i(\phi_1 - \phi_0)} \| \\
    = 0
\end{split}
\end{equation}

\section{Proof 2}
\label{sec:sup_proof_2}
The similarity between two vectors of phases encoded into the complex state of an oscillator can be computed via superposition without directly recovering the relative phase between them. 

This arises from the fact that the half-angle of the relative phase between the two complex potentials is trigonometrically related to the hypotenuse of the addition (superposition) of the two complex states. 

\begin{equation}
\begin{split}
    cos(\frac{1}{2}(\phi_1 - \phi_0)) = \frac{1}{2} \| Z_1 + Z_0 \| \\
    \therefore (\phi_1 - \phi_0) = 2 \cdot arccos(\frac{\| Z_1 + Z_0 \|}{2}) \\
    cos(\phi_1 - \phi_0) = cos(2 \cdot arccos[\frac{\| Z_1 + Z_0 \|}{2}])
\end{split}
\end{equation}

\section{Proof 3}
\label{sec:sup_proof_3}
The binding operation produces an angle that is the sum of two input angles. When these angles are known, this operation is trivial. 

\begin{equation}
    \phi_2 = \phi_0 + \phi_1
\end{equation}

However, when these angles are encoded into the complex-valued state of an oscillator, this operation becomes more involved. Angular information can be encoded in the temporally invariant differences between two oscillators’ instantaneous phase. In order to accurately decode these angles, a ``reference oscillator'' must be used to define the ``starting point'' for the angle expressed by the complex state of an oscillator. 

\begin{equation}
\begin{split}
    \phi_{ref} = 0 \\
    Z_{ref} = e^{i(\omega t + \phi_{ref})} = e^{i \omega t} \\
    Z_0 = e^{i(\omega t + \phi_0)} \\
    Z_1 = e^{i(\omega t + \phi_1)} \\
    Z_2 = e^{i(\omega t + \phi_0 + \phi_1)} \\
\end{split}
\end{equation}

The desired product of binding two oscillators is the complex state $Z_2$. A geometric construction shows that this point can be reached by following two chords along the complex unit circle determined by $\phi_0$ and $\phi_1$. This ``journey'' begins at the reference oscillator’s state and is displaced by the difference between it and the first oscillator’s state. 

\begin{equation}
    Z' = Z_{ref} + (Z_0 - Z_{ref}) = Z_0
\end{equation}

The second step to reach the final state is defined by the difference $Z_1-Z_{ref}$. However, this step must be ``rotated'' by the value of $\phi_0$ to begin at $Z_0$ and end at the desired point $Z_2$. This angle is already encoded in $Z_0$; multiplying by this value in the complex domain achieves the required rotation. 

\begin{equation}
    Z'' = Z_0 \cdot (Z_1 - Z_{ref})
\end{equation}

When the representations for $Z$ are expanded, however, one can  see that the multiplication of two time-varying complex values together causes the frequency at which they are rotating to double. 

\begin{equation}
\begin{split}
    Z'' = e^{i(\omega t + \phi_0)} \cdot (e^{i(\omega t + \phi_1) - e^{i(\omega t + \phi_0)}}) \\
    = e^{i(2\omega t + \phi_0 + \phi_1)} -  e^{i(2\omega t + \phi_0 )}\\
\end{split}
\end{equation}

This frequency doubling can be counteracted by adding a third multiplicative term consisting of the reference’s complex conjugate.

\begin{equation}
\begin{split}
    Z''' = e^{i(\omega t + \phi_0)} \cdot (e^{i(\omega t + \phi_1) - e^{i(\omega t + \phi_0)}}) \cdot e^{-i\omega t}\\
    = e^{i(\omega t + \phi_0 + \phi_1)} -  e^{i(\omega t + \phi_0 )}
\end{split}
\end{equation}

Thus, the values $Z'$ and $Z'''$ can  be summed to follow the path of two chords rotating at identical rates to reach the output of the binding operation on complex oscillators.

\begin{equation}
\begin{split}
    Z_2 = Z' + Z''' \\
    = Z_0 + Z_0 \cdot (Z_1 - Z_{ref}) \cdot \overline{Z_{ref}}
\end{split}
\end{equation}

While similarity and bundling calculations natively depend on differences of phases that remain invariant w.r.t. time, binding requires decoding the original phase encoded into the oscillators and thus requires a reference.

\section{Neural Network Input Transformations}
\label{sec:nn_transforms}

An array of 21 by 13 silicon pixels is simulated in the ``Smart Pixels'' dataset. Each pixel records the amount of charge accumulated within its volume as a pion passes through the sensor ($Q(x,y,t)$). The charge level is recorded every 5 picoseconds for 20 steps, giving a recording time of 100 ps. 

Independent splits of the data are used to calibrate and normalize the amount of charge that is induced within the pixels. The 99.99th percentile of these charges is approximately 15,000 kiloelectrons (ke). This is taken as the 6$\sigma$ value; all charge values are normalized by $6 / 15,000$, and values above this threshold are clipped to stabilize the training process.

\begin{equation}
    Q'(x,y,t) = \frac{6 \cdot Q(x,y,5)}{15,000 ke}
\end{equation}

Indirect-drive neural networks utilize the normalized charge $Q'$ at the end of the sampling period (100 ps) as the input to the neural network. Charges across each column are summed to produce the ``y-profile'' of charge produced by the particle.

\begin{equation}
    X_{Q, indirect} = sum(Q'(i,1:21, 100 \space ps)), i \in 1:13
\end{equation}

In the direct-drive neural network, the difference in normalized charge levels between each time step is taken as approximating the current value during that period.

\begin{equation}
    I(x,y,t) = \frac{Q'(x,y,t+5 ps) - Q'(x,y,t)}{5 \space ps}
\end{equation}

Currents through the entire 100 ps recording period are used as inputs that excite oscillators. Each input directly activates a corresponding neuron. Additionally, a small Gaussian ``timing pulse'' is added halfway through the integration period at 50 ps to ensure that each neuron will produce an output to drive the downstream network.

\begin{equation}
    \begin{split}
    GK(\mu, \sigma, t) = exp(-1 \cdot \frac{t - \mu}{2 \cdot \sigma}^{2}) \\
    I_{timing}(t) = GK(\mu = 5 ps, \sigma = 0.01 \space ps, t) \\
    I_{Q, direct}(x,y,t) = sum(I(i,1:21,t)) + I_{timing}(t), i \in 1:13
    \end{split}
\end{equation}

The span of the simulated pixel detector is small with respect to the range in curvatures the simulated particles take. To allow each particle to ``hit'' the simulated detector, it is shifted to meet the point at which the particle meets the detector ``plane.'' This shift in the y direction provides crucial information on the particle's momentum and must be transmitted as an input into the network to provide an efficacious classification. The pixel sensor can be ``shifted'' within $\pm 32.5$ mm. The y location is thus normalized by this value to produce a phase.

\begin{equation}
    \begin{split}
        Y'(y) = y / 32.5 \space mm \\
        X_{Y, direct}(y) = Y'(y) / 2 + 0.5
    \end{split}
\end{equation}

Similarly to the clock pulse, this y-location phase can be converted into to a short current pulse for the direct-drive network.

\begin{equation}
    I_{Y, direct}(y, t, T) = GK(\mu = X_{Y, direct}(y), \sigma = 0.01 \space ps, t \% T)
\end{equation}
\noindent
Here $T$ is the period of the oscillator being driven and $\%$ is the modulo operation. Together with the y-profile information the input for each classification becomes a vector of 14 values in the case of the indirect-drive network and a vector of 14 time-varying functions in the case of the direct-drive network. 

\end{appendices}

\bibliographystyle{IEEEtran}
\bibliography{IEEEabrv,HDC.bib}

\end{document}